\documentclass[12pt,]{amsart}
\usepackage{amsmath}
\usepackage{amscd}
\usepackage{amssymb}
\usepackage{latexsym}
\input xy
\xyoption{all} \CompileMatrices

\newtheorem{theorem}{Theorem}[section]
\newtheorem{lemma}[theorem]{Lemma}
\newtheorem{corollary}[theorem]{Corollary}
\newtheorem{proposition}[theorem]{Proposition}

\newtheorem{remark}[theorem]{Remark}
\newtheorem{definition}[theorem]{Definition}

\newcommand{\bgl}{\begin{equation}}         
\newcommand{\egl}{\end{equation}}
\newcommand{\bgln}{\begin{eqnarray}}        
\newcommand{\egln}{\end{eqnarray}}
\newcommand{\bglnoz}{\begin{eqnarray*}}     
\newcommand{\eglnoz}{\end{eqnarray*}}
\newcommand{\btheo}{\begin{theorem}}
\newcommand{\etheo}{\end{theorem}}
\newcommand{\blemma}{\begin{lemma}}
\newcommand{\elemma}{\end{lemma}}
\newcommand{\bproof}{\begin{proof}}
\newcommand{\eproof}{\end{proof}}
\newcommand{\bbew}{\begin{beweis}}
\newcommand{\ebew}{\end{beweis}}
\newcommand{\bremark}{\begin{remark}\em}
\newcommand{\eremark}{\end{remark}}
\newcommand{\bdefin}{\begin{definition}}
\newcommand{\edefin}{\end{definition}}
\newcommand{\bprop}{\begin{proposition}}
\newcommand{\eprop}{\end{proposition}}
\newcommand{\bcor}{\begin{corollary}}
\newcommand{\ecor}{\end{corollary}}
\newcommand{\mn}{\par\medskip\noindent}

\newcommand{\cB}{\mathcal B}
\newcommand{\cC}{\mathcal C}
\newcommand{\cD}{\mathcal D}

\newcommand{\cJ}{\mathcal J}
\newcommand{\cK}{\mathcal K}
\newcommand{\cL}{\mathcal L}

\newcommand{\cO}{\mathcal O}

\newcommand{\cQ}{\mathcal Q}

\newcommand{\lori}{\longrightarrow}

\newcommand{\ve}{\varepsilon}
\newcommand{\vp}{\varphi}

\def\SEMI{\mbox{$\times\kern-2pt\vrule height5pt width.6pt \kern3pt $}}

\newcommand{\Ker}{{\rm Ker\,}}
\newcommand{\Aut}{{\rm Aut\,}}

\newcommand{\id}{{\rm id}}

\def\Az{\mathbb{A}}
\def\Cz{\mathbb{C}}
\def\Fz{\mathbb{F}}

\def\Nz{\mathbb{N}}

\def\Qz{\mathbb{Q}}
\def\Rz{\mathbb{R}}
\def\Zz{\mathbb{Z}}

\def\Tz{\mathbb{T}}
\def\A{\mathfrak{A}}
\def\Af{\mathfrak{A}[\varphi]}
\def\AR{\mathfrak{A}[R]}

\begin{document}

\title{C*-algebras associated with algebraic actions}
\date{\today}
\author[J. Cuntz]{Joachim Cuntz}
\address{Mathematisches Institut, Einsteinstr.62, 48149
M\"unster, Germany} \email{cuntz@uni-muenster.de}
\thanks{Research supported by DFG through
CRC 878 and by ERC through AdG 267079}
\subjclass[2000]{Primary: 22D25, 46L80, 46L89, 11R04, 11M55, 37A55} \keywords{endomorphism, compact abelian group,
C*-algebra, purely infinite, $K$-theory, algebraic integers, semigroup, KMS-states}

\begin{abstract}\noindent This is a survey of work in which the author was involved in recent years. We consider C*-algebras constructed from representations of one or several algebraic endomorphisms of a compact abelian group - or, dually, of a discrete abelian group. In our survey we do not try to describe the entire scope of the methods and results obtained in the original papers, but we concentrate on the important thread coming from the action of the multiplicative semigroup of a Dedekind ring on its additive group. Representations of such actions give rise to particularly intriguing problems and the study of the corresponding C*-algebras has motivated many of the new methods and general results obtained in this area.\end{abstract}\maketitle

\section{Introduction}
By an algebraic action we mean here an action of a semigroup by algebraic endomorphisms of a compact abelian group or, dually, by endomorphisms of a discrete abelian group. Such actions are much studied in ergodic theory but they also give rise to interesting C*-algebras. In fact, quite a few of the standard examples of simple C*-algebras such as $\cO_n$-algebras, Bunce-Deddens algebras, UHF-algebras etc. arise from canonical representations of such endomorphisms. But the class of C*-algebras obtained from general algebraic actions is much vaster and exhibits new interesting phenomena.

We start our survey with the discussion, following \cite{CuVe}, of the C*-algebra $\mathfrak A[\alpha]$ generated by the so called Koopman representation on $L^2(H)$ of a single endomorphism $\alpha$, satisfying natural conditions, of a compact abelian group $H$, together with the natural representation of $C(H)$. This C*-algebras is always simple purely infinite and can be described by a natural set of generators and relations. It contains a canonical maximal commutative C*-algebra $\cD$ with spectrum a Cantor space. This subalgebra is generated by the range projections $s^ns^{n*}$, where $s$ is the isometry implementing the given endomorphism, and by their conjugates $u_\gamma s^ns^{n*}u_\gamma^*$, under the unitaries $u_\gamma$ given by the characters $\gamma$ of $H$. Then, the subalgebra $\cB$ generated by $\cD$ together with the $u_\gamma$ is of Bunce-Deddens type and simple with unique trace. Moreover, $\mathfrak A[\alpha]$ can be considered as a crossed product of $\cB$ by a single endomorphism.

The next case we consider is the C*-algebra generated analogously by the Koopman representation of a family of commuting endomorphisms. We consider the important special case of endomorphisms arising from the ring of integers $R$ in a number field $K$. The multiplicative semigroup $R^\times$ acts by commuting endomorphisms on the additive group $R \cong \Zz^n$ or equivalently on the dual group $\widehat{R} \cong \Tz^n$ ($n$ being the degree of the field extension $K$ over $\Qz$). The commutative semigroup $R^\times$ has a non-trivial structure and acts by interesting endomorphisms on $\Tz^n$. The study of the C*-algebra $\AR$ generated by the Koopman representation in this situation goes back to \cite{Cun} and was originally motivated by connections to Bost-Connes systems \cite{BoCo}.

Again, $\AR$ is simple purely infinite and is described by natural generators and relations. It has analogous subalgebras $\cD$ and $\cB$, and $\AR$ can be viewed as a semigroup crossed product $\cB \rtimes R^\times$. The new and challenging problem is the computation of the $K$-theory of $\AR$. The key to this computation is a duality result for adele-groups and corresponding crossed products, \cite{CuLi2}.

Since $\AR$ is generated by the Koopman representation, on $\ell^2R$, of the semidirect product semigroup $R\rtimes R^\times$, the next very natural step in our program is the consideration of the C*-algebra generated by the natural representation of this semigroup on $\ell^2(R \rtimes R^\times)$ rather than on $\ell^2R$, i.e. of the left regular C*-algebra $C^*_\lambda(R \rtimes R^\times)$. This algebra is still purely infinite but no longer simple. It can be described by natural generators and relations. The algebra $\AR$ is a quotient of $C^*_\lambda(R \rtimes R^\times)$ and the latter algebra is defined by relaxing the relations defining $\AR$ in a systematic way. The best way to do so is to add a family of projections, indexed by the ideals of the ring $R$, as additional generators and to incorporate those into the relations. This way of defining the relations also guided Xin Li in his description of the left regular C*-algebras for more general semigroups \cite{LiSG}.

The (non-trivial) problem of computing the $K$-theory of $C^*_\lambda(R \rtimes R^\times)$ turned out to be particularly fruitful \cite{CEL1}, \cite{CEL2}. It led to a powerful new method for computing the $K$-groups, for regular C*-algebras of more general semigroups and of crossed products by automorphic actions of such more general semigroups, as well as for crossed products of certain actions of groups on totally disconnected spaces. In the special case of $C^*_\lambda(R \rtimes R^\times)$ we get the interesting result that the $K$-theory is described by a formula that involves the basic number theoretic structure of the number field $K$, namely the ideal class group and the action of the group of units (invertible elements in $R$) on the additive group of an ideal.

Finally, we include a brief discussion of the rich $KMS$-structure on $C^*_\lambda(R \rtimes R^\times)$ for the natural one-parameter action on this C*-algebra. Just as the $K$-theory for $C^*_\lambda(R \rtimes R^\times)$, this structure is related to the number theoretic invariants of $R$, resp. $K$.

Our goal in this survey is limited. We try to describe a leitmotif in this line of research and to explain the connections and similarities between the various results. The original articles contain much more information and many additional finer, more sophisticated and more general results which we omit completely. We also do not describe the results in the order they were obtained originally, but rather in the order which seems more systematic with hindsight.

\section{Single algebraic endomorphisms}\label{sec:1}

Let $H$ be a compact abelian group and $G=\widehat{H}$ its dual discrete
group. We assume that $G$ is countable. Let
$\alpha$ be a surjective endomorphism of $H$ with finite kernel. We
denote by $\vp$ the dual endomorphism $\chi\mapsto \chi\circ \alpha$ of
$G$ (i.e. $\vp=\widehat{\alpha}$). By duality, $\vp$ is injective and
has finite cokernel, i.e. the quotient $G/\vp G$ will be finite.
Both $\alpha$ and $\vp$ induce isometric endomorphisms $s_\alpha$
and $s_\vp$ of the Hilbert spaces $L^2H$ and $\ell^2G$,
respectively. This isometric representation of $\alpha$ on $L^2(H)$ is called the Koopman representation in ergodic theory.

We will also assume that $$\bigcap_{n\in\Nz}\vp^nG=\{0\}$$ which, by
duality, means that $$\bigcup_{n\in\Nz}\Ker \alpha^n$$ is dense in
$H$ (this implies in particular that $H$ and $G$ can not be finite). These conditions on $\alpha$ are quite natural and for instance apply to the usual examples considered in ergodic theory. We list a few important examples of compact groups and endomorphisms satisfying our conditions at the end of this section.

We want to describe the C*-algebra
$C^*(s_\alpha,C(H))$ generated in $\cL(L^2H)$ by $C(H)$, acting by
multiplication operators, and by the isometry $s_\alpha$. Via
Fourier transform it is isomorphic to the C*-algebra
$C^*(s_\vp,C^*G)$ generated in $\cL(\ell^2 G)$ by $C^*G$, acting via
the left regular representation, and by the isometry $s_\vp$. These
two unitarily equivalent representations are useful for different
purposes.

Now, $C^*(s_\vp,C^*G)$ is generated by the isometry $s=s_\vp$ together
with the unitary operators $u_g,\, g\in G$ and these operators satisfy the relations

\bgl\label{rel}u_gu_h=u_{g+h}\qquad s u_g=u_{\vp (g)}s\qquad
\sum_{g\in G/\vp G}u_gss^*u_g^*=1\egl

\bdefin Let $H,G$ and $\alpha, \vp$ be as above. We denote by $\Af$
the universal C*-algebra generated by an isometry $s$ and unitary
operators $u_g,\, g\in G$ satisfying the relations
(\ref{rel}).\edefin

It is shown in \cite{CuVe} that $\Af\cong C^*(s_\alpha,C(H))\cong
C^*(s_\vp,C^*G)$, i.e. that the natural map from the universal C*-algebra to the C*-algebra generated by the concrete Koopman representation is an isomorphism. Particular situations of interest arise when $H = (\Zz/n)^\infty$ with $\alpha$ the left shift (this gives rise to $\Af \cong \cO_n$) or when $H = \Tz^n$.

\blemma\label{D} The C*-subalgebra $\cD$ of $\Af$ generated by all
projections of the form $u_gs^ns^{*n}u_g^*$, $g\in G, n\in \Nz$ is
commutative. Its spectrum is the ``$\vp$-adic completion''
$$G_\vp=\mathop{\lim}\limits_{ {\scriptstyle\longleftarrow_n} }G/\vp^nG$$
It is an inverse limit of the finite spaces $G/\vp^nG$ and becomes a
Cantor space with the natural topology.

$G$ acts on $\cD$ via $d\mapsto u_gdu_g^*$, $g\in G,\, d\in \cD$.
This action corresponds to the natural action of the dense subgroup
$G$ on its completion $G_\vp$ via translation. The map $\cD\to\cD$
given by $x\mapsto sxs^*$ corresponds to the map induced by $\vp$ on
$G_\vp$.\elemma

From now on we will denote the compact abelian group $G_\vp$ by $M$. By construction, $G$ is a dense subgroup of $M$.
The dual group of $M$ is the discrete abelian group
$$L = \mathop{\lim}\limits_{ {\scriptstyle\longrightarrow_n}}\Ker
(\alpha^n :H\to H)$$ Because of the condition that we impose on
$\alpha$, $L$ can be considered as a dense subgroup of $H$.

The groups $M$ and $L$ play an important role in the analysis
of $\Af$. They are in a sense complementary to $H$ and $G$. By Lemma
\ref{D}, the C*-algebra $\cD$ is isomorphic to $C(M)$ and to
$C^*(L)$.

\btheo\label{B} The C*-subalgebra $B_\vp$ of $\Af$ generated by
$C(H)$ together with $C(M)$ (or equivalently by $C^*G$ together with
$C^*L$) is isomorphic to the crossed product $C(M)\rtimes G$. It is
simple and has a unique trace. \etheo
\btheo\label{pure} The algebra $\Af$ is simple, nuclear and purely
infinite. Moreover, it is isomorphic to the semigroup crossed product
$B_\vp\rtimes_{\gamma_\vp} \Nz$ (i.e. to the universal unital
C*-algebra generated by $B_\vp$ together with an isometry $t$ such
that $txt^*=\gamma_\vp (x)$, $x\in B_\vp$).\etheo
The fact that $\Af$ is a crossed product $B_\vp\rtimes \Nz$ can be used to prove the following

\btheo\label{KA} (cf. \cite{CuVe}) The $K$-groups of $\Af$ fit into an exact sequence
as follows
\bgl\label{PVO}\xymatrix{K_*C(H)\ar[r]^{1-b(\vp)}&
K_*C(H)\ar[r]&K_*\Af\ar@/^5mm/[ll]}\egl\mn
where the map $b( \vp):K_*C(H) \to K_*C(H)$ satisfies $b( \vp) \alpha_* = N(\alpha)\id$ with $N(\alpha) := |\Ker \alpha |$.
\etheo
In \cite{CuVe}, the analysis of $\Af$ and the formula \eqref{PVO} for its $K$-theory was also extended to the case where $\alpha$ is replaced by a so called rational polymorphism.

There are quite a few papers in the literature containing special cases or parts of the results described in this section. We mention only \cite{Hirsh} where it was shown that $\Af$ is simple and characterized by generators and relations and \cite{EHR} where in particular a formula similar to \eqref{PVO} was derived for an expansive endomorphism of $\Tz^n$ - both papers using methods different from \cite{CuVe}.
\subsection{Examples} Here are some examples of endomorphisms in the class we consider. \mn

1. Let $H=\prod_{k\in \Nz}\Zz/n$,
$G=\bigoplus_{k\in \Nz}\Zz/n$ and $\alpha$ the one-sided shift on
$H$ defined by $\alpha ((a_k))=(a_{k+1})$.

We obtain $M=\prod_{k\in \Nz}\Zz/n\cong H$ and $L=\bigoplus_{k\in
\Nz}\Zz/n\cong G$. The algebra $B_\vp$ is a UHF-algebra of type
$n^\infty$ and $\Af$ is isomorphic to $\cO_n$. It is interesting to
note that the UHF-algebra $B_\vp$ is generated by two maximal
abelian subalgebras both isomorphic to $C(M)$.

2. Let $H=\Tz$, $G=\Zz$ and $\alpha$ the endomorphism of
$H$ defined by $\alpha (z)=z^n$. The algebra $B_\vp$ is a
Bunce-Deddens-algebra of type $n^\infty$ and $\Af$ is isomorphic to
a natural subalgebra of the algebra $\cQ_\Nz$ considered in
\cite{Cun}. In this case, we also get for $B_\vp$ the interesting
isomorphism $C(\Zz_n)\rtimes \Zz\cong C(\Tz)\rtimes L$ where $\Zz$
acts on the Cantor space $\Zz_n$ by the odometer action (addition of
1) and $L$ denotes the subgroup of $\Tz$ given by all $n^k$-th roots
of unity, acting on $\Tz$ by translation.

3. Let $H=\Tz^n$, $G=\Zz^n$ and $\alpha$ an endomorphism
of $H$ determined by an integral matrix $T$ with non-zero
determinant. We assume that the condition
$$\bigcap_{n\in\Nz}\vp^nG=\{0\}$$ is satisfied (this is in fact not very restrictive).

The algebra $B_\vp$ is a higher-dimensional analogue of a
Bunce-Deddens-algebra. In the case where $H$ is the additive group
of the ring $R$ of algebraic integers in a number field of degree
$n$ and the matrix $T$ corresponds to an element of $R$, the algebra
$\Af$ is isomorphic to a natural subalgebra of the algebra $\A [R]$
considered in the following section. It is also isomorphic to the algebra
studied in \cite{EHR}.

4. As another natural example related to number theory consider the additive group of the polynomial ring $\Fz_p[t]$ over a finite field. An endomorphism satisfying our conditions is given by multiplication by a non-zero element in $\Fz_p[t]$. In this case $\Af$ is related to certain graph C*-algebras, see \cite{CuLiP}.

5. Let  $p$ and $q$ be natural numbers that
are relatively prime and $\gamma$ the endomorphism of $\Tz$ defined
by $z\mapsto z^p$. We take

$$H=\mathop{\lim}\limits_{ \mathop{{\scriptstyle\longleftarrow}}
\limits_\gamma}\Tz\quad\qquad G=\Zz[\frac {1}{p}]$$\mn
$\alpha_q$ the endomorphism of $H$ induced by $z\mapsto z^q$ and
$\vp_q$ the endomorphism of $G$ defined by $\vp_q(x)=qx$. These
endomorphisms satisfy our hypotheses. We find that $M=\Zz_q$ (the $q$-adic
completion of $\Zz$).\mn

In all these examples one can work out the $K$-theory of $\Af$ using formula \eqref{PVO}, see \cite{CuVe}.
\section{Actions by a family of endomorphisms, ring C*-algebras}\label{sec:2}
It is a natural problem to extend the results of section \ref{sec:1} to actions of a family (semigroup) of several commuting endomorphisms of a compact abelian group,  satisfying the conditions of section \ref{sec:1}. It turns out that the structural results such as simplicity, pure infiniteness, canonical subalgebras carry over without problem. However the computation of the $K$-groups needs completely new ideas.

The most prominent example for us arises as follows. Let $K$ be a number field, i.e. a finite algebraic extension of $\Qz$. The ring of algebraic integers $R \subset K$ is defined as the integral closure of $\Zz$ in $K$, i.e. as the set of elements $a \in K$ that annihilate some monic polynomial with coefficients in $\Zz$. This ring is always a Dedekind domain (a Dedekind domain is by definition an integral domain in which every nonzero proper ideal factors into a product of prime ideals). It has many properties similar to the ordinary ring of integers $\Zz \subset \Qz$, but it is not a principal ideal domain in general. Its additive group is always isomorphic to $ \Zz^n$ where $n$ is the degree of the field extension.

Consider the multiplicative semigroup $R^\times = R\setminus \{0\}$ of $R$. It acts as endomorphisms on the additive group $R$ and thus also on the compact abelian dual group $\widehat{R}\cong \Tz^n$. Such endomorphisms of $\Tz^n$ are a frequent object of study in ergodic theory. If $R$ is not a principal ideal domain, the semigroup $R^\times$ has an interesting structure.

As in section \ref{sec:1} we consider the Koopman representation of $R^\times$ on $L^2(\widehat{R}) \cong \ell^2R$.

\bdefin We denote by $\AR$ the C*-algebra generated by $C(\widehat{R})$ and $R^\times$ on $L^2(\widehat{R})$ (or equivalently the C*-algebra generated by the action of $C^*(R)$ and of $R^\times$ on $\ell^2R$). \edefin

$\AR$ is generated by the isometries $s_n$, $n \in R^\times$ and the unitaries $u_j$, $j \in R$. The $s_n$ define a representation of the abelian semigroup $R^\times$ by isometries, the $u_j$ define a representation of the abelian group $R$ by unitaries and together they satisfy the relations

\bgl\label{AR} s_nu_k=u_{kn}s_n,\;k \in R, \, n,m \in
R^\times\qquad \sum_{j \in R/nR}u_js_ns_n^*u_{-j}\;=1  \egl

The basic analysis of the structure of $\AR$ is completely parallel to the discussion in section \ref{sec:1} (in fact historically the article \cite{CuLi} preceded \cite{CuVe}). One obtains
\btheo (cf. \cite{CuLi}) The C*-algebra $\AR$ is simple purely infinite and nuclear. It is the universal C*-algebra generated by a unitary representation $u$ of $R$ together with an isometric representation $s$ of $R^\times$ satisfying the relations \eqref{AR}. \etheo
As for $\Af$ in section \ref{sec:1} there are canonical subalgebras $\cD$ and $\cB$ of $\AR$. The spectrum of the  commutative C*-algebra $\cD$ is a Cantor space canonically homeomorphic to the maximal compact subring of the space of finite adeles for the number field $K$. The subalgebra $\cB$ is generated by $\cD$ together with the $u_j$, $j \in R$. It is simple and has a unique trace (a higher dimensional Bunce-Deddens type algebra). The general structure of C*-algebras associated like this with a ring has been developed further by Xin Li in \cite{LiRing}.

In order to compute the $K$-groups for $\AR$ the natural strategy would appear to be an iteration of the formula \eqref{PVO} of Theorem \ref{KA}. Since the proof of formula \eqref{PVO} is based on the usual Pimsner-Voiculescu sequence this would amount to iterating this sequence in order to compute the $K$-groups for the crossed product by $\Zz^n$ by a commuting family of $n$ automorphisms. However this strategy immediately runs into problems since, assuming the $K$-groups for the crossed product by the first automorphism are determined, it is not at all clear how the second automorphism will act on these groups. In other words, there is a spectral sequence abutting to the $K$-theory for the crossed product by $\Zz^n$, but it is useless for actual computations without further knowledge of the higher boundary maps in the spectral sequence. An analysis of relevant properties of the spectral sequence for actions as here is contained in \cite{Bar}.

The key to the computation of the $K$-groups for $\AR$ in \cite{CuLi2} is the following duality result.
\btheo\label{du} Let $\Az_f$ and $\Az_\infty$ denote the locally compact spaces of finite, resp. infinite adeles of $K$ both with the natural action of the additive group $K$. Then the crossed product C*-algebras $C_0(\Az_f) \rtimes K$ and $C_0(\Az_\infty) \rtimes K$ are Morita equivalent, equivariantly for the action of $K^\times$ on both algebras. \etheo

Note that the space $\Az_\infty$ is simply $\Rz^n$ where $n$ is the degree of the field extension. From this theorem the $K$-groups of $\AR$ can be computed, at least in the case where the only roots of unit in $K$ are $\pm 1$.

We explain this here only for the case where $K= \Qz$, $R = \Zz$. In this case everything becomes rather concrete. The spectrum of the canonical commutative subalgebra $\cD$ is the profinite completion $\overline{\Zz}$ of $\Zz$ (we use here $\overline{\Zz}$ rather than the more standard notation $\widehat{\Zz}$ in order not to create confusion with the dual group of $\Zz$). It is homeomorphic to the infinite product of the $p$-adic completions $\Zz_p$ for all primes $p$ in $\Zz$. Moreover $\Az_f$ is the restricted infinite product of the $\Qz_p$ and $\Az_\infty$ simply is $\Rz$

Thus Theorem \ref{du} gives a Morita equivalence between $C_0(\Az_f)\rtimes \Qz$ and $C_0(\Rz)\rtimes \Qz$. Moreover, by a standard argument, the first crossed product is Morita equivalent to $\cB \cong C(\overline{\Zz})\rtimes \Zz$.

Denote by $\cB'$ the C*-algebra generated by $\cB$ together with the symmetry $s_{-1}$, i.e. $\cB' \cong \cB \rtimes \Zz/2$ for the action of $s_{-1}$. Since $\cB' \cong (C(\overline{\Zz})\rtimes \Zz)\rtimes \Zz/2$ is an inductive limit of $C((\Zz/n\Zz)\rtimes \Zz)\rtimes \Zz/2$ and this latter algebra is isomorphic to $M_n(C^*(\Zz \rtimes \Zz/2))$ it is not difficult to compute the $K$-theory of $\cB$ as $K_0(\cB') = \Zz \oplus \Qz$ and $K_1(\cB')=0$.

Now, we can use the Pimsner-Voiculescu sequence to compute the $K$-theory of the crossed product $\mathfrak A_1 = \cB' \rtimes \Nz = C^*(\cB', s_2)$ as
$$K_0(\mathfrak A_1) = \Zz \qquad K_1(\mathfrak A_1) = \Zz$$

By a slight refinement of the statement in Theorem \ref{du}, $\mathfrak A_1$ is Morita equivalent to $(C_0(\Rz)\rtimes \Qz)\rtimes (\Zz/2 \times \Zz)$ where $\Zz/2 \times \Zz$ acts by multiplication by $-1$ and by $2$.

Denote now by $\mathfrak A_n$ the C*-algebra generated by $\cB'$ together with the $s_{p_1},\ldots, s_{p_n}$, where $p_1, \ldots , p_n$ denote the first $n$ prime numbers (with $p_1=2$). Then again $\mathfrak A_n$ is Morita equivalent to $(C_0(\Rz)\rtimes \Qz)\rtimes (\Zz/2\times\Zz^n)$, where $\Zz/2$ acts by multiplication by $-1$ and $\Zz^n$ by multiplication by $p_1, \ldots , p_n$. Moreover $\AR$ is the inductive limit of the $\mathfrak A_n$.

We can now consider the canonical inclusions
\bgl\label{iota}\iota_n:C_0(\Rz)\rtimes (\Zz/2\times\Zz^n) \to (C_0(\Rz)\rtimes \Qz)\rtimes (\Zz/2\times\Zz^n)\, \sim_{\textrm{Morita}}\;\mathfrak A_n\egl
into the crossed product where we leave out the action of the additive $\Qz$ by translation on the left hand side.

By the discussion above, $\iota_1$ induces an isomorphism in $K$-theory. Now we obtain $\iota_{n+1}$ from $\iota_n$ by taking the crossed product by $\Zz$ (acting by multiplication by $p_{n+1}$) on both sides in \eqref{iota}. Therefore, applying the Pimsner-Voiculescu sequence on both sides, we deduce, using the five-lemma, from the fact that $\iota_n$ induces an isomorphism on $K$-theory that the same holds for $\iota_{n+1}$.  The important point is that the action of $\Zz^n$ on the left hand side is homotopic to the trivial action, simply because multiplication by $p_1, \ldots , p_n$ is homotopic to multiplication by 1 on $\Rz$. Therefore $K_*( \mathfrak A_n) \cong K_*((C_0(\Rz)\rtimes\Zz/2) \otimes C^*\Zz^n)$.

As a consequence we obtain
\btheo (\cite{CuLi2} The map $\iota_n$ induces an isomorphism on K-theory for all $n$. The $K$-theory of $\AR$ is isomorphic to the $K$-theory of $(C_0(\Rz)\rtimes \Zz/2)\rtimes \Qz^\times$. \etheo
Note that the $K$-theory of $C_0(\Rz)\rtimes \Zz/2$ is the same as the one of $\Cz$ and that therefore the $K$-theory of $\AR$ is the same as the one of an infinite-dimensional torus.

The argument that we sketched for $K=\Qz$ works in a very similar, though somewhat more involved way for a number field with $\pm 1$ as only roots of unit. In this case one has to determine the $K$-theory of $C_0(\Rz^n)\rtimes \Zz/2$ rather than that of $C_0(\Rz)\rtimes \Zz/2$.
The case of an arbitrary number field $K$ can be treated in the same fashion. The important difference comes from the more general group $\mu (K)$ of roots of unit. For the computation one needs non-trivial information on the $K$-theory of the crossed product $C_0(\Rz^n) \rtimes \mu(K)$ and thus on the equivariant $K$-theory of $\Rz^n$ with respect to the action of $\mu(K)$. This non-trivial computation has been carried through by Li and L\"uck in \cite{LiLu} using previous work by Langer and L\"uck \cite{LaLu}.

The analysis of the structure and of the $K$-theory of $\AR$ can also be carried out in the case where $R$ is a polynomial ring over a finite field (ring of integers in a certain function field). The structure of the C*-algebra in this case is more closely related to the example of the shift endomorphism of $(\Zz/p\Zz)^\infty$ mentioned above and to certain Cuntz-Krieger algebras. Nevertheless for the computation of the $K$-theory one can again use the duality result in Theorem \ref{du} and the result for the $K$-theory is again similar, \cite{CuLiP}.

\section{Regular C*-algebra for `ax+b'-semigroups}

By definition, the ring C*-algebra $\AR$ discussed in section \ref{sec:2} is obtained from the natural representations of $C^*(R) \cong C(\widehat{R})$ and of the semigroup $R^\times$ on the Hilbert space $\ell^2 R \cong L^2(\widehat{R})$. Another way to view this is to say that it is defined by the natural representation of the semidirect product semigroup $R \rtimes R^\times$ on $\ell^2R$.

Now, this semidirect product semigroup has an even more natural representation, given by the left regular representation on the Hilbert space $\ell^2 (R \rtimes R^\times )$. The study of the left regular C*-algebra $C^*_\lambda(R \rtimes R^\times )$ was begun in \cite{CDL}. This C*-algebra is no longer simple but still purely infinite and has an intriguing structure. In particular, it has a very interesting $KMS$-structure and the determination of its $K$-theory leads to new challenging problems.

The first obvious observation concerning $C^*_\lambda(R \rtimes R^\times )$ is that, just as $\AR$, it is generated by a unitary representation $u_x, x\in R$ of the additive group $R$ and a representation by isometries $s_a, a\in R^\times$ of the multiplicative semigroup $R^\times$ satisfying the additional relation $s_au_x = u_{ax}s_a$. However the last relation $ \sum_{x \in R/aR}u_xs_as_a^*u_{-x}=1 $ in \eqref{AR} becomes
\bgl\label{reg}  \sum_{x \in R/aR}u_xs_as_a^*u_{-x} \;\leq 1  \egl
In fact, it turns out that this weakened relation \eqref{reg} (of course together with the relations on the $u_x, s_a$ in the previous paragraph) determines $C^*_\lambda(R \rtimes R^\times)$ in the case where $R$ is a principal ideal domain. The general case however is more intricate. In general, it is still possible to describe $C^*_\lambda(R \rtimes R^\times)$ by natural defining relations. However the most natural way to do so uses an incorporation of the natural idempotents obtained as range projections of the partial isometries given by products of the $u_x$, $s_a$ and their adjoints. It turns out that these range projections correspond exactly to the ideals of $R$.

The generators singled out in \cite{CDL} then are
$u_x, x\in R$; $s_a,\, a \in
R^\times$; $e_I$, $I$ a non-zero ideal in $R$. The relations are
\begin{enumerate}
  \item The $u_x$ are unitary and satisfy $u_xu_y=u_{x+y}$, the
  $s_a$ are isometries and satisfy $s_as_b=s_{ab}$. Moreover $s_a u_x=u_{ax}s_a$ for all $x\in R,\,a\in
  R^\times$.
  \item The $e_I$ are projections and satisfy $e_{I\cap J}=e_Ie_J$, $e_R=1$.
  \item We have $s_ae_I = e_{aI}s_a$.
  \item For $x\in I$ one has $u_xe_I = e_Iu_x$, for $x\notin I$
  one has $e_Iu_xe_I=0$.
\end{enumerate}
The universal C*-algebra with these generators and relations is no longer simple but in most respects its structure is similar to the one of $\AR$. There is a canonical maximal commutative subalgebra $\cD$ with totally disconnected spectrum (generated by the $e_I$), and a Bunce-Deddens type subalgebra $\cB$ generated by $\cD$ together with the $u_x$, $x \in R$. Using this structure one shows
\btheo\label{Clam} (cf. \cite{CDL}) The universal C*-algebra with generators $u_x$, $s_a$, $e_I$ satisfying the relations 1.,2.,3.,4. above is canonically isomorphic to $C^*_\lambda(R \rtimes R^\times)$. As a consequence $C^*_\lambda(R \rtimes R^\times)$ is also isomorphic to the semigroup crossed product $\cD \rtimes (R \rtimes R^\times)$ (i.e. to the universal C*-algebra generated by $\cD$ together with a representation of the semigroup $R\rtimes R^\times$ by isometries implementing the given endomorphisms of $\cD$).\etheo
It follows that $\AR$ is a quotient of $C^*_\lambda(R \rtimes R^\times )$. As mentioned above, in the simplest case $\AR$ is obtained from $C^*_\lambda(R \rtimes R^\times )$ by `tightening' the relation $\sum_{j \in R/nR}u_js_ns_n^*u_{-j}$ $\leq 1$ to $\sum_{j \in R/nR}u_js_ns_n^*u_{-j} \;= 1$. This kind of tightening has occurred in many places in the literature under the name tight representation or boundary quotient etc.

The relations 1.,2.,3.,4. above  turned out to also give the right framework for describing the left regular C*-algebra of more general semigroups. The theory of these regular C*-algebras has been developed by Xin Li \cite{LiSG}, \cite{LiNuc}.

As in section \ref{sec:2} the key to the computation of $K_*(C^*_\lambda(R \rtimes R^\times ))$, for the ring $R$ of integers in a number field $K$, lies in a $KK$-equivalence between the given action by endomorphisms of our semigroup with a basically trivial situation.

The semigroup $S=R \rtimes R^\times $ admits $G=K \rtimes K^\times$ as a canonical enveloping group. The action of $S$ on the commutative subalgebra $\cD$ of $C^*_\lambda(R \rtimes R^\times )$ has a natural dilation to an action of $G$. This means that $\cD$ can be embedded into a larger commutative C*-algebra $\cC \supset \cD$ with an action of $G$ which extends the action of $S$ on $\cD$ (this uses the fact that $S$ is a directed set ordered by right divisibility). The crossed product $\cC \rtimes G$ is then Morita equivalent to $\cD \rtimes S \cong C^*_\lambda(R \rtimes R^\times)$ (the last isomorphism follows from Theorem \ref{Clam}).

A fractional ideal in $K$ is a subset $J$ of $K$ of the form $J=aI$ where $I$ is an ideal in $R$ and $a \in K^\times$. Denote by $\cJ$ the set of all fractional ideals of $R$ in $K$, i.e. the set of all translates in $K$ of ideals in $R$ under the action of $K^\times$.

It is easy to see that, in the dilated system, there is a bijection $J \mapsto e_J$ between $\cJ$ and the translates under $G$ of the projections $e_I$, $I$ ideal in $R$. Moreover the $e_J$, $J \in \cJ$ generate $\cC \supset \cD$.
Using the fact that $R$ is a Dedekind domain it is not difficult to show that the family $\{e_J\}$ forms a regular basis of $\cC$ in the sense of the following definition. The importance of the regularity condition (or, in another guise, of the `independence' of the family of constructible left ideals of the semigroup) has been noticed by Xin Li.

\begin{definition}\label{reg1} If $\{e_J:J\in \cJ\}$ is a countable set of non-zero projections in a commutative C*-algebra $\cC$, we say that $\{e_J\}$ is a {\em regular basis}
for $\cC$ if it is linearly independent, closed under multiplication (up to $0$) and generates $\cC$ as a C*-algebra (this means that  ${\rm span\,}\{e_J:J\in \cJ\}$ is a dense subalgebra of $\cC$).
\end{definition}
Now the group $G$ acts on $\cC$, on $\cJ$ and on the algebra $\cK=\cK(\ell^2(\cJ))$ of compact operators. We can trivially define an equivariant *-homomorphism $\kappa : C_0(\cJ) \to \cK\otimes \cC$ by mapping $\delta_J$ to $\ve_I \otimes e_J$. Here, $\delta_J$ denotes the indicator function of the one-point set $\{J\}$ and $\ve_J$ denotes the matrix in $\cK$ which is is 1 in the diagonal place $(J,J)$ and 0 otherwise (matrix unit).

\begin{theorem}\label{kap}(\cite{CEL1,CEL2})
The equivariant map $\kappa$ induces an isomorphism
$$K_*(C_0(\cJ)\rtimes G) \lori K_*(\cC\rtimes G) \cong K_*(C^*_\lambda(R \rtimes R^\times))$$
\end{theorem}

But now, by Green's imprimitivity theorem the crossed product $C_0(\cJ)\rtimes G$ is simply Morita equivalent to the direct sum, over the $G$-orbits in
$\cJ$, of the C*algebras of the stabilizer groups of each orbit.
\bdefin\label{class} The ideal class group $Cl_K$ is the quotient of the semigroup of fractional ideals in $K$ under the equivalence relation where $J$ is equivalent to $J'$ iff there is $a \in K^\times$ such that $J'=aJ$. \edefin
If $R$ is the ring of algebraic integers in the number field $K$, then the class group is a finite abelian group.

Two fractional ideals $J$ and $J'$ are in the same orbit for $G$ if and only if there is $a \in K^\times$ such that $J'= aJ$. Therefore, by Definition \ref{class} the orbits are labeled  exactly by the elements of the class group $Cl_K$. The stabilizer group of the class of a fractional ideal $J$ then is given by the semidirect product $J \rtimes R^*$ of the additive group $J$ by the group $R^*$ of units (i.e. of invertible elements in $R^\times$). As a corollary to Theorem \ref{kap} we thus obtain
\bcor\label{Kred} For each element $\gamma$ of the class group $Cl_K$ choose any ideal $I_\gamma$ representing the class $\gamma$. Then
$$ K_*(C^*_\lambda(R \rtimes R^\times))\cong \bigoplus_{\gamma\in  Cl_{K}}K_*\big(C^*(I_\gamma)\rtimes R^*\big).$$
\ecor

In the situation at hand, Theorem \ref{kap} can be proven directly - essentially in a similar way as at the end of section \ref{sec:2} using the equivariant map $\kappa$. There is however a much more powerful approach based on techniques from work on the Baum-Connes conjecture developed by Echterhoff and others \cite{ENO},\cite{CEO}. They establish the following principle:
\begin{itemize}
\item[]
Assume that the group $G$ satisfies the Baum-Connes conjecture with
coefficients in the $G$-algebras $A$ and $B$. Let $\kappa : A \to B$ be an equivariant homomorphism which induces, via descent, isomorphisms $K_*(A\rtimes
H)\cong K_*(B\rtimes H)$ for all compact subgroups $H$ of $G$. Then
$\kappa$ also induces an isomorphism $K_*(A\rtimes_rG)\cong
K_*(B\rtimes_rG)$.
\end{itemize}
Theorem \ref{kap} then follows from checking that the equivariant map $C_0(\cJ) \to \cK \otimes \cC$ used there satisfies this condition for all finite subgroups of $G$.

This approach to Theorem \ref{kap} has a much broader scope of applications. It allows to extend the argument to general actions of a group $G$, that satisfies the Baum-Connes conjecture with coefficients, on a commutative C*-algebra $\cC$ admitting a $G$-invariant regular basis of projections in the sense of Definition \ref{reg1}. In particular, it can then be used to compute the $K$-theory of the left regular C*-algebra for a large class of semigroups as well as for crossed products by automorphic actions by such semigroups. Moreover this more general method also allows to compute the $K$-theory for crossed products for an action of a group on a totally disconnected space that admits a regular basis as in Definition \ref{reg1}, \cite{CEL1},\cite{CEL2}.

For instance, one obtains
\begin{theorem}\label{thm-Dedekind}(\cite{CEL2})
 Let $R$ be a Dedekind domain with quotient field $Q(R)$ and $A$ a C*-algebra. Then the following are true:
 \begin{enumerate}
 \item For every action $\alpha:R^\times \to \Aut(A)$ there is a canonical isomorphism
 $$ K_*(A\rtimes_{\alpha,r} R^\times)\cong \bigoplus_{\gamma\in Cl_{Q(R)}}K_*(A\rtimes_{\alpha, r} R^*).$$
 \item For every action $\alpha: R^\times/R^*\to \Aut(A)$  there is a canonical isomorphism
 $$ K_*(A\rtimes_{\alpha,r} (R^\times/R^*))\cong \bigoplus_{\gamma\in Cl_{Q(R)}}K_*(A).$$
 \item For every action $\alpha:R\rtimes R^\times\to \Aut(A)$ there is a canonical isomorphism
 $$ K_*\big(A\rtimes_{\alpha,r}(R\rtimes R^\times)\big)\cong \bigoplus_{\gamma\in  Cl_{Q(R)}}K_*\big(A\rtimes_{\alpha,r}(I_\gamma\rtimes R^*)\big).$$
 \end{enumerate}
 \end{theorem}

 The above method of computing $K$-theory for semigroup C*-algebras and for certain crossed products for actions on totally disconnected spaces has been developed further by Li-Norling in \cite{LiNor}, \cite{LiNor2}.

\section{$KMS$-states}
To end this survey we briefly discuss the $KMS$-structure for the natural one-parameter automorphism group of $C^*_\lambda (R \rtimes R^\times )$ where, again, $R$ is the ring of algebraic integers in a number field $K$. After all, part of the motivation for the study of ring C*-algebras came from Bost-Connes system and a main feature of such systems is the rich $KMS$-structure. Also, one of the reasons in \cite{CDL} for passing from the ring C*-algebra $\AR$ to $C^*_\lambda (R \rtimes R^\times )$ was the existence of many $KMS$-states on the latter algebra.

Recall that, for a non-zero ideal $I$ in $R$, we denote by $N(I)$ the
norm of $I$, i.e. the number $N(I)=|R/I|$ of elements in $R/I$. For
$a\in R^\times$ we also write $N(a)=N(aR)$. The norm is
multiplicative, \cite{Neukirch}. Using the norm one defines a natural one-parameter automorphism
group $(\sigma_t)_{t\in\Rz}$ on $C^*_\lambda (R \rtimes R^\times )$, given on the generators by
$$\sigma_t(u_x)=u_x\quad\sigma_t(e_I)=e_I\quad\sigma_t(s_a)=N(a)^{it}s_a$$
(this assignment manifestly respects the relations between the
generators and thus induces an automorphism). Let $\beta$ be a real number $\geq 0$. Recall that a
$\beta$-$KMS$ state with respect to a one parameter automorphism group
$(\sigma_t)_{t\in\Rz}$ is a state $\vp$ which satisfies $\vp
(yx)=\vp (x\sigma_{i\beta}(y))$ for a dense set of analytic vectors
$x,y$ and for the natural extension of $(\sigma_t)$ to complex
parameters on analytic vectors, \cite{BrRo}. For the
one-parameter automorphism group $\sigma$, defined above, the $\beta$-$KMS$ condition for a state $\vp$
translates to \bgl\label{KMScond} \vp (u_xz)=\vp (z\,u_x)\quad\vp
(e_I z)=\vp (z\,e_I)\quad\vp (s_az)=N(a)^{-\beta}\vp (z\,s_a)\egl for a
set of analytic vectors $z$ with dense linear span and for the
standard generators $u_x$, $e_I$, $s_a$ of $C^*_\lambda (R \rtimes R^\times )$.

\btheo\label{kms} (\cite{CDL}) The KMS-states on $C^*_\lambda (R \rtimes R^\times )$ at inverse temperature $ \beta$ can be described. One has

\begin{enumerate}
  \item no KMS-states for $ \beta < 1$.
  \item for each $ \beta \in [1,2 ]$ a unique $ \beta$-KMS state.
  \item for $ \beta \in (2, \infty)$ a bijection between $ \beta$-KMS states and traces on

      $$ \bigoplus_{ \gamma \in Cl_K} C^*(I_\gamma) \rtimes R^* $$

      where $Cl_K$ is the ideal class group, $I_\gamma$ is any ideal representing $ \gamma$ and $R^*$ denotes the multiplicative group of invertible elements in $R$ (units).
\end{enumerate}
\etheo
The simpler case of Theorem \ref{kms}, where $R=\Zz, \,K=\Qz$, had essentially been treated already by Laca-Raeburn in \cite{LaRa}. The first assertion in Theorem \ref{kms} is basically obvious. The proof, in \cite{CDL}, of point 3. uses special representations of $C^*_\lambda (R \rtimes R^\times )$ which seem to be of independent interest. The proof of 2. in \cite{CDL} uses a result, also of some independent interest, on asymptotics of partial Dedekind $\zeta$-functions. An alternative subsequent proof of Theorem \ref{kms}, due to Neshveyev, is obtained by relating the problem to a general result on $KMS$-states for C*-algebras of non-principal groupoids, \cite{Nesh}.

There is a striking parallel between the formula for the $KMS$-states for $\beta > 2$ in the theorem above and the formula for the $K$-theory of $C^*_\lambda (R \rtimes R^\times )$ in Corollary \ref{Kred}. The $K$-theory is isomorphic to the $K$-theory of the C*-algebra $ \bigoplus_{ \gamma \in \Gamma} C^*(I_\gamma) \rtimes R^* $ while the simplex of $KMS$-states is in bijection with the traces of this direct sum C*-algebra. Note that both results are quite non-trivial, as $ \bigoplus_{ \gamma \in \Gamma} C^*(I_\gamma) \rtimes R^* $ is not a natural subalgebra of $C^*_\lambda (R \rtimes R^\times )$.

\end{document}